\pdfoutput=1
\RequirePackage{ifpdf}
\ifpdf 
\documentclass[pdftex]{sigma}
\else
\documentclass{sigma}
\fi

\numberwithin{equation}{section}
\newtheorem{Theorem}{Theorem}[section]
\newtheorem{Lemma}[Theorem]{Lemma}
{\theoremstyle{definition}
\newtheorem{Definition}[Theorem]{Definition}
\newtheorem{Remark}[Theorem]{Remark}
}

\begin{document}

\allowdisplaybreaks

\newcommand{\arXivNumber}{1409.4287}

\renewcommand{\PaperNumber}{116}

\FirstPageHeading

\ShortArticleName{Representation Theory for Conf\/luent Cherednik Algebras}

\ArticleName{Non-Symmetric Basic Hypergeometric Polynomials\\
and Representation Theory for Conf\/luent\\
Cherednik Algebras}

\Author{Marta MAZZOCCO}

\AuthorNameForHeading{M.~Mazzocco}

\Address{Department of Mathematical Sciences, Loughborough University, Loughborough LE11 3TU, UK}
\Email{\href{mailto:m.mazzocco@lboro.ac.uk}{m.mazzocco@lboro.ac.uk}}
\URLaddress{\url{http://homepages.lboro.ac.uk/~mamm4/}}

\ArticleDates{Received October 31, 2014, in f\/inal form December 19, 2014; Published online December 30, 2014}

\Abstract{In this paper we introduce a~basic representation for the conf\/luent Cherednik algebras $\mathcal H_{\rm V}$,
$\mathcal H_{\rm III}$, $\mathcal H_{\rm III}^{D_7}$ and $\mathcal H_{\rm III}^{D_8}$ def\/ined in~arXiv:1307.6140.
To prove faithfulness of this basic representation, we introduce the non-symmetric versions of the continuous
dual $q$-Hahn, Al-Salam--Chihara, continuous big $q$-Hermite and continuous $q$-Hermite polynomials.}


\Classification{33D80; 33D52; 16T99}

\section{Introduction}

In this paper we introduce a~faithful representation on the space $\mathcal A$ of Laurent polynomials for the conf\/luent
Cherednik algebras $\mathcal H_{\rm V}$, $\mathcal H_{\rm III}$, $\mathcal H_{\rm III}^{D_7}$ and $\mathcal H_{\rm III}^{D_8}$
def\/ined in~\cite{M}\footnote{See \cite[Theorem~4.1]{M} for $\mathcal H_{\rm V}$, $\mathcal H_{\rm III}$ and \cite[Def\/inition~1.4]{M} for
$\mathcal H_{\rm III}^{D_7}$ and $\mathcal H_{\rm III}^{D_8}$~--
observe that in~\cite{M}~$W$ is~$X^{-1}$ for these
algebras.} as conf\/luences of the Cherednik algebra of type $\check{C_1}C_1$~\cite{Cher,NS,Sa}:
\begin{itemize}\itemsep=0pt
\item $\mathcal H_{\rm V}$ is the algebra generated by $T_0$, $T_1$, $X^{\pm 1}$ with relations:
\begin{gather}
\label{sahi2-V}
(T_1+ab)(T_1+1)=0,
\\
\label{sahi3-V}
T_0(T_0+1)=0,
\\
\label{sahi4-V}
(T_1 X +a)(T_1 X +b)=0,
\\
\label{sahi6-V}
q T_0 X^{-1}+ c = X (T_0+1).
\end{gather}
\item $\mathcal H_{\rm III}$ is the algebra generated by $T_0$, $T_1$, $X^{\pm 1}$ with relations:
\begin{gather}
\label{sahi2-III}
(T_1+ab)(T_1+1)=0,
\\
\label{sahi3-III}
T_0^2=0,
\\
\label{sahi4-III}
(T_1 X +a)(T_1 X +b)=0,
\\
\label{sahi6-III}
q T_0 X^{-1}+1= X T_0.
\end{gather}
\item $\mathcal H_{\rm III}^{D_7}$ is the algebra generated by $T_0$, $T_1$, $X^{\pm 1}$ with relations:
\begin{gather}
\label{sahiD71}
T_1(T_1+1)=0,
\\
\label{sahiD72}
T_0^2=0,
\\
\label{sahiD73}
T_1 X +a-X^{-1}(T_1+1) =0,
\\
\label{sahiD74}
q T_0 X^{-1}+ 1-X T_0 =0.
\end{gather}
\item $\mathcal H_{\rm III}^{D_8}$ is the algebra generated by $T_0$, $T_1$, $X^{\pm 1}$ with relations:
\begin{gather}
\label{sahiPIIID82}
T_1(T_1+1)=0,
\\
\label{sahiPIIID83}
T_0^2=0,
\\
\label{sahiPIIID84}
T_1 X -X^{-1} (T_1+1) =0,
\\
\label{sahiPIIID85}
qT_0 X^{-1}+ 1-X T_0=0.
\end{gather}
\end{itemize}

To prove faithfulness of our basic representation (see Theorems~\ref{th:repPV},~\ref{th:repPIII} and~\ref{th:repD7} here
below) in each case, we select a~special basis of polynomials in $\mathcal A$ on which the operators (or specif\/ic
combinations of them) act nicely.
These bases are obtained by considering the non-symmetric versions of the continuous dual $q$-Hahn, Al-Salam--Chihara,
continuous big $q$-Hermite and continuous $q$-Hermite polynomials respectively.

In~\cite{Sa} Sahi introduced the non-symmetric version of Koornwinder polynomials~\cite{AW}, and proved that they form
a~basis in the space $\mathcal A$ of Laurent polynomials.
A~detailed discussion of the rank one case, i.e.~the non-symmetric Askey--Wilson polynomials, was presented in~\cite{NS}
(see also~\cite{K1}).
It turns out that these non-symmetric Askey--Wilson polynomials behave well under the subsequent degeneration limits
$d\to 0$, $c\to 0$, $b\to 0$ and f\/inally $a\to 0$.
However the proof of faithfulness of our basic representation is not a~straightforward limit of the same proof in the
case of the Askey--Wilson algebra, as one would naively expect.
This is because the f\/irst degeneration limit destroys some of the leading coef\/f\/icients in the positive powers of~$z$ of
half the non-symmetric continuous dual $q$-Hahn polynomials and their degenerations.
Moreover, the algebra $\mathcal H_{\rm III}$ is not in fact the limit of $\mathcal H_{\rm V}$ as $c\to 0$ but the one as
$c\to\infty$, which introduces the need of an isomorphism and a~few tricks.
Last but not least, the $\mathcal H_{\rm III}^{D_7}$ and $\mathcal H_{\rm III}^{D_8}$ do not admit a~presentation {\it \`a~la
Bernstein--Zelevinsky}, which makes the proof of faithfulness in that case rather involved.

\section[Non-symmetric continuous dual $q$-Hahn polynomials and basic representation of $\mathcal H_{\rm V}$]{Non-symmetric
continuous dual $\boldsymbol{q}$-Hahn polynomials\\ and basic representation of $\boldsymbol{\mathcal H_{\rm V}}$}
\label{se:PV}

The continuous dual $q$-Hahn polynomials are the following (we write them here in monic form like in~\cite{KLS}):
\begin{gather*}
p_n(z;a,b,c):= \frac{(ab,ac;q)_n}{a^n} \, {}_3\phi_2\left(
\begin{matrix}
q^{-n},az,a z^{-1}
\\
ab,ac
\end{matrix}
;q,q \right),
\end{gather*}
and can be obtained from the Askey--Wilson polynomials as limits when $d\to 0$.
This same limit can be performed on the non-symmetric the Askey--Wilson polynomials, leading to the following (here we
follow~\cite{K1} approach):
\begin{Definition}
Let
\begin{gather*}
q_n^\dagger (z;a,b,c):= q^{\frac{n-1}{2}} (z-c) p_{n-1}\big(q^{-\frac{1}{2}}z;q^{\frac{1}{2}}a,q^{\frac{1}{2}}b,q^{\frac{1}{2}} c\big),
\end{gather*}
the non-symmetric continuous dual $q$-Hahn polynomials are def\/ined as follows:
\begin{gather*}
E_{-n}[z]:= p_n(z;a,b,c)-q_n^\dagger (z;a,b,c),
\qquad
n=1,2,\dots,
\\
E_{n}[z]:=q^n p_n(z;a,b,c)+\big(1-q^n\big) q_n^\dagger (z;a,b,c),
\qquad
n=1,2,\dots,
\\
E_{0}[z]:= 1.
\end{gather*}
\end{Definition}

\begin{Theorem}
\label{th:repPV}
For $q,a,b,c\neq 0$, $q^m\neq 1$ $(m=1,2,\dots)$, the algebra $\mathcal H_{\rm V}$ has a~faithful representation on the space
$\mathcal A$ of Laurent polynomials $f[z]$ as follows:
\begin{gather}
\label{t0PV}
(T_0 f)[z]:=\frac{(z-c)z}{q-z^2} \left(f[z] - f\left[q z^{-1}\right]\right),
\\
\label{t1PV}
(T_1 f)[z]:=\frac{(a+b)z-(1+a b)}{1-z^2} f[z] + \frac{(1-a z)(1-b z)}{1-z^2} f\left[z^{-1}\right],
\\
\label{XPV}
(X f)[z]:=z f[z].
\end{gather}
\end{Theorem}

To prove this theorem we follow the same outline as the proof of Theorem~5.3 in~\cite{K1} with some important changes as
explained in Remark~\ref{rmk:difference1} here below.

First of all, to prove that the operators def\/ined by~\eqref{t0PV}--\eqref{XPV} satisfy the relations~\eqref{sahi2-V}--\eqref{sahi6-V} is a~straightforward computation.
To prove faithfulness, we need the following two lemmata:

\begin{Lemma}
Let
\begin{gather}
\label{eq:ZY}
Z:=(T_0+1)T_1^{-1}
\qquad
\hbox{and}
\qquad
Y:=T_1 T_0,
\end{gather}
the algebra $\mathcal H_{\rm V}$ can equivalently be described as the algebra generated by $T_1$, $X^{\pm 1}$, $Y$, $Z$, satisfying
the following relations respectively:
\begin{gather}
\label{LD00-PV}
Z Y=Y Z=0,
\\
\label{LD1-PV}
X T_1 = -a bT_1^{-1} X^{-1} -a-b,
\\
\label{LD2-PV}
T_1^{-1} Y = Z T_1 - 1,
\\
\label{LD3-PV}
(T_1+ab)(T_1+1)=0,
\\
\label{LD4-PV}
ab Y X =- q T_1^2 X Y - q(a+b) T_1 Y -a b T_1 X +ab c T_1.
\end{gather}
The algebra $\mathcal H_{V}$ is spanned by elements $X^m Y^n T_1^i$ and $X^m Z^n T_1^i$, where $m\in\mathbb Z$,
$n\in\mathbb N$ and $i=1,2$.
\end{Lemma}

\begin{proof}
To prove the equivalence it is enough to observe that by def\/ining~$Z$ and~$Y$ as in~\eqref{eq:ZY},
relations~\eqref{LD00-PV}--\eqref{LD4-PV} follow
from~\eqref{sahi2-V}--\eqref{sahi6-V}.
Vice-versa, def\/ining $T_0:= T_1^{-1}Y$ we see that relations~\eqref{LD00-PV}--\eqref{LD4-PV} imply~\eqref{sahi2-V}--\eqref{sahi6-V}.

To prove that $\mathcal H_{V}$ is spanned by elements $X^m Y^n T_1^i$ and $X^m Z^n T_1^i$, where $m\in\mathbb Z$,
$n\in\mathbb N$ and $i=1,2$, we use the relations~\eqref{LD00-PV}--\eqref{LD4-PV} and the further relations which can be obtained as a~consequence
of~\eqref{LD00-PV}--\eqref{LD3-PV}:
\begin{gather*}
Y X^{-1} = q^{-1} X^{-1} Y + q^{-1}(1+ ab) X^{-1} Z T_1 - q^{-1}(a+b) Z T_1
+q^{-1} X^{-1} T_1-c q^{-1} T_1,
\\
Z X=q^{-1} X Z - q\frac{1+a b}{a b} X^{-1} Z T_1+\frac{a+ b}{a b} Z T_1- \frac{1}{a b} X^{-1} T_1
\nonumber
\\
\phantom{Z X=}{}
+\frac{c}{a b q} T_1+ \frac{(1+a b)(q-1)}{a b} \left(X^{-1}-\frac{c}{q}\right)
\end{gather*}
to order any word in the algebra as wanted.
\end{proof}

\begin{Lemma}
The non-symmetric continuous dual $q$-Hahn polynomials form a~basis in the space $\mathcal A$ of Laurent polynomials and
are eigenfunctions of the operators $Y:= T_1T_0$ and $Z:=(T_0+1)T_1^{-1}$:
\begin{gather}
\label{eq:eigenY}
(Y E_{-n})[z]= \frac{1}{q^n} E_{-n}[z],
\qquad
n=1,2,3,\dots,
\\
(Y E_{n})[z]=0,
\qquad
n=0,1,2,\dots.
\nonumber
\\
\label{eq:eigenZ}
(Z E_{-n})[z]=0,
\qquad
n=1,2,3,\dots,
\\
(Z E_{n})[z]=- \frac{1}{a b q^n} E_{n}[z],
\qquad
n=0,1,2,\dots.
\nonumber
\end{gather}
\end{Lemma}

\begin{proof}
By using the def\/inition of the $q$-hypergeometric series $ {}_3\phi_2$ it is easy to prove that the terms with the
highest powers in~$z$ and $\frac{1}{z}$ in $E_{-n}$ and $E_n$ have the following form
\begin{gather}
\label{andamentim}
E_{-n}[z]= z^{-n}+\dots+ \big(a b c q^{n-1}-a-b\big) z^{n-1},
\qquad
n=1,2,\dots,
\\
\label{andamentip}
E_{n}[z]=z^n+\dots + q^n z^{-n},
\qquad
n=1,2,\dots.
\end{gather}
Using these relation it is straightforward to prove that the non-symmetric continuous dual $q$-Hahn polynomials form
a~basis in $\mathcal A$.

Now to prove~\eqref{eq:eigenY}, we use the fact that the operator~$Y$ acts on $\mathcal A$ as follows
\begin{gather*}
(Y f)[z]:=\frac{(z-c)z(1-(a+ b)z+a b)}{(1-z^2)(q-z^2)}\big(f\big[q z^{-1}\big]-f[z]\big)
\nonumber
\\
\phantom{(Y f)[z]:=}{}
+ \frac{(1-a z)(1-b z)(1-c z)}{(1-z^2)(1-q z^2)}\big(f[q z]- f\big[z^{-1}\big]\big).
\end{gather*}
Observe that thanks to the forward shift operator relation (14.3.8) in~\cite{KLS}, one has:
\begin{gather*}
q_n^\dagger (z;a,b,c)= -\frac{q^n z(z-c)}{(q^n-1)(q-z^2)} \big(p_n(z;a,b,c)-p_n \big(q^{-1} z;a,b,c\big)\big),
\end{gather*}
so that one can express $(Y E_{-n})[z]- \frac{1}{q^n} E_{-n}[z]$ only in terms of $p_n(z;a,b,c)$, $p_n (q z;a,b,c)$ and
$p_n (q^{-1} z;a,b,c)$, which can be shown to be zero by using the $q$-dif\/ference equation (14.3.7) in~\cite{KLS}.
In a~similar manner all other relations are proved.
\end{proof}

\begin{Remark}
\label{rmk:difference1}
Note that as shown in~\eqref{andamentim}, the polynomials $E_{-n}[z]$ do not have a~term of order $z^{n}$ like the
non-symmetric Askey--Wilson polynomials did.
This is due to the fact that the coef\/f\/icient of the term $z^n$ in the non-symmetric Askey--Wilson polynomials tends to
zero as $d\to 0$.
The absence of such term makes the end of the proof of Theorem~\ref{th:repPV} more tricky than proof of Theorem~5.3  in~\cite{K1}.
\end{Remark}

\begin{proof}[Proof of Theorem~\ref{th:repPV}.] First by using the symmetry properties of the continuous dual $q$-Hahn polynomials
and their properties it is easy to show that
\begin{gather*}
(T_1 E_{-j}) [z]=- \big(1 +ab - ab q^j\big)E_{-j}[z] - ab E_j[z],
\\
(T_1 E_{j}) [z]=\big(1-q^j\big)\big(1-a b q^j\big)E_{-j}[z] - ab q^{j}E_j[z].
\end{gather*}
Combining this with~\eqref{eq:eigenY}--\eqref{andamentip}, we can prove the
following $\forall\, n>0$, $\forall\, m\in\mathbb Z$, $\forall\, j>0$:
\begin{gather}
X^m E_{-j}[z]= z^{m- j}+\dots+ \big(a b c q^{j-1}-a-b\big) z^{m+j-1},
\nonumber
\\
X^m Y^n E_{-j}[z]= q^{-j n} z^{m- j}+\dots+ q^{-j n} \big(a b c q^{j-1}-a-b\big) z^{m+j-1},
\nonumber
\\
X^m Y^n T_1 E_{-j}[z]=- \big(1 +ab - ab q^j\big)q^{-j n}\big(z^{m-j}+\dots + \big(a b c q^{j-1}-a-b\big) z^{m+j-1}\big),
\nonumber
\\
X^m T_1 E_{-j}[z]= - (1 +ab) z^{m-j}+\dots - ab z^{m+j},
\nonumber
\\
X^m Z^n E_{-j}[z]= 0,
\nonumber
\\
X^m Z^n T_1 E_{-j}[z]=\left(\frac{-1}{a b}\right)^{n-1} q^{-n j} \big(z^{m+j}+\dots + q^j z^{m-j}\big),
\nonumber
\\
X^m E_{j}[z]= z^{m+j}+\dots+ q^{j}z^{m-j},
\label{eq-ej-z}
\\
X^m Y^n E_{j}[z]= 0,
\nonumber
\\
X^m Y^n T_1 E_{j}[z]=\big(1-q^j\big) \big(1 - ab q^j\big)q^{-j n}\big(z^{m-j}+\dots + \big(a b c q^{j-1}-a-b\big) z^{m+j-1}\big),
\nonumber
\\
X^m T_1 E_{j}[z]= \big(1 - abq^j-q^j\big) z^{m-j}-a b q^j z^{m+j},
\nonumber
\\
X^m Z^n E_{j}[z]=\left(\frac{-1}{a b q^j}\right)^{n} \big(q^j z^{m-j}+\dots + z^{m+j}\big),
\nonumber
\\
X^m Z^n T_1 E_{j}[z]=\left(\frac{-1}{a b q^j}\right)^{n-1} \big(z^{m+j}+\dots + q^j z^{m-j}\big).
\nonumber
\end{gather}
Now assume by contradiction that a~linear combination acts as zero operator in our representation, let us write such
linear combination as:
\begin{gather*}
\sum\limits_m a_m X^m+ \sum\limits_{m,n} b_{m,n} X^m Y^n+\sum\limits_{m,n, i}c_{m,n} X^m Y^n T_1+\sum\limits_{m,n} d_{m}X^m T_1
\\
\qquad{}+\sum\limits_{m,n}e_{m,n} X^m Z^n +\sum\limits_{m,n,}f_{m,n} X^m Z^n T_1.
\end{gather*}
Take the minimum value~$M$ of~$m$ such that at least one coef\/f\/icient $a_{m}$, $b_{m,n}$, $c_{m,n}$, $d_{m}$, $e_{m,n}$, $f_{m,n}$
is nonzero.
Acting on $E_j$, and collecting the terms with the minimum possible power of~$z$, by~\eqref{eq-ej-z} we obtain the
equation:
\begin{gather*}
a_M q^j+\sum\limits_n c_{M,n}\big(1-q^j\big)\big(1-a b q^j\big)q^{-jn}+d_M \big(1-a b q^j-q^j\big)
\\
\qquad
{}+\sum\limits_n e_{M,n}\left(\frac{-1}{a b q^j}\right)^{n} q^j+\sum\limits_n f_{M,n}\left(\frac{-1}{a b q^j}\right)^{n-1}
q^j=0,
\qquad
\forall\, j>0.
\end{gather*}
It is easy to prove that for generic values of the parameters $a$, $b$, $c$, this is an inf\/inite set of linearly independent
equations, therefore the only possible solution is the trivial one.
So we can only have coef\/f\/icients of type $b_{M,n}$ not zero.
Again, acting on $E_{-j}$, and collecting the terms with the minimum possible power of~$z$ we obtain for every $j>0$,
the equation:
\begin{gather*}
\sum\limits_n b_{M,n}q^{-jn}=0,
\end{gather*}
which admit only trivial solutions.
\end{proof}

\section[Non-symmetric Al-Salam--Chihara polynomials and basic representation of $\mathcal H_{\rm III}$]{Non-symmetric
Al-Salam--Chihara polynomials\\ and basic representation of $\boldsymbol{\mathcal H_{\rm III}}$}

The Al-Salam--Chihara polynomials are the following:
\begin{gather*}
Q_n(z;a,b):= \frac{(ab;q)_n}{a^n}\, {}_3\phi_2\left(
\begin{matrix}
q^{-n},az,a z^{-1}
\\
ab,0
\end{matrix}
;q,q \right),
\end{gather*}
and can be obtained from the continuous dual $q$-Hahn polynomials as limits when $c\to 0$.
By taking the limit $c\to 0$ of the non-symmetric continuous dual $q$-Hahn polynomials we obtain the following:

\begin{Definition}
Let
\begin{gather*}
Q_n^\dagger (z;a,b):= q^{\frac{n-1}{2}} z   Q_{n-1}\big(q^{-\frac{1}{2}} z;q^{\frac{1}{2}}a,q^{\frac{1}{2}}b\big),
\end{gather*}
the non-symmetric Al-Salam--Chihara polynomials are def\/ined as follows:
\begin{gather*}
E_{-n}[z]:= Q_n(z;a,b)-Q_n^\dagger (z;a,b),
\qquad
n=1,2,\dots,
\\
E_{n}[z]:=q^n Q_n(z;a,b)+\big(1-q^n\big) Q_n^\dagger (z;a,b),
\qquad
n=1,2,\dots,
\\
E_{0}[z]:= 1.
\end{gather*}
\end{Definition}

\begin{Theorem}
\label{th:repPIII}
For $q,a,b\neq 0$, $q^m\neq 1$ $(m=1,2,\dots)$, the algebra $\mathcal H_{\rm III}$ has a~faithful representation on the
space $\mathcal A$ of Laurent polynomials $f[z]$ as follows:
\begin{gather}
\label{t0PIII}
(T_0 f)[z]:=-\frac{z}{q-z^2} \left(f[z] - f\left[q z^{-1}\right]\right),
\\
\label{t1PIII}
(T_1 f)[z]:=\frac{(a+b)z-(1+a b)}{1-z^2} f[z] + \frac{(1-a z)(1-b z)}{1-z^2} f\left[z^{-1}\right],
\\
\label{XPIII}
(X f)[z]:=z f[z].
\end{gather}
\end{Theorem}

To prove that the operators def\/ined by~\eqref{t0PIII}--\eqref{XPIII} satisfy the
relations~\eqref{sahi2-III}--\eqref{sahi6-III} is a~straightforward computation.
To prove faithfulness, we again need to provide an equivalent representation for the algebra $\mathcal H_{\rm III}$.
This is where we need to be careful as the relation between the non-symmetric Al-Salam--Chihara polynomials and the
algebra $\mathcal H_{\rm III}$ is not as straightforward as before because the algebra $\mathcal H_{\rm III}$ was
obtained as limit of $\mathcal H_{\rm V}$ as $c\to\infty$ rather than $c\to 0$.
However, changing the def\/inition of~$Z$ and~$Y$ we can still prove the following:

\begin{Lemma}
Let
\begin{gather}
\label{eq:new-YZ}
Z:= - X T_0T_1^{-1}+ T_1^{-1}
\qquad
\hbox{and}
\qquad
Y:= -T_1 X T_0,
\end{gather}
then the algebra $\mathcal H_{\rm III}$ can equivalently be described as the algebra generated by $T_1$, $X^{\pm 1}$, $Y$, $Z$,
satisfying the following relations respectively:
\begin{gather}
\label{LD00-PIII}
Z Y=Y Z=0,
\\
\label{LD1-PIII}
X T_1 = -a bT_1^{-1} X^{-1} -a-b,
\\
\label{LD2-PIII}
T_1^{-1} Y = Z T_1 - 1,
\\
\label{LD3-PIII}
(T_1+ab)(T_1+1)=0,
\\
\label{LD4-PIII}
ab Y X =- q T_1^2 X Y - q(a+b) T_1 Y -a b T_1 X.
\end{gather}
The algebra $\mathcal H_{\rm III}$ is spanned by elements $X^m Y^n T_1^i$ and $X^m Z^n T_1^i$, where $m\in\mathbb Z$, $n\in\mathbb N$ and $i=1,2$.
\end{Lemma}

\begin{proof}
The relations~\eqref{LD00-PIII}--\eqref{LD4-PIII} follow
from~\eqref{sahi2-III}--\eqref{sahi6-III}.
Vice-versa, def\/ining $T_0:= - X^{-1} T_1^{-1}Y$ we see that
relations~\eqref{LD00-PIII}--\eqref{LD4-PIII}
imply~\eqref{sahi2-III}--\eqref{sahi6-III}.

To prove that $\mathcal H_{\rm III}$ is spanned by elements $X^m Y^n T_1^i$ and $X^m Z^n T_1^i$, where $m\in\mathbb Z$,
$n\in\mathbb N$ and $i=1,2$, we use the relations~\eqref{LD00-PIII}--\eqref{LD4-PIII} and the further equivalent relations
\begin{gather*}
Y X^{-1} = q^{-1} X^{-1} Y + q^{-1}(1+ab) X^{-1} Z T_1 - q^{-1}(a+b) Z T_1
+q^{-1} X^{-1} T_1,
\\
Z X=q^{-1} X Z - q\frac{1+a b}{a b} X^{-1} Z T_1+\frac{1+a b}{a b} Z T_1- \frac{1}{a b} X^{-1} T_1
+ \frac{(1+a b)(q-1)}{a b} X^{-1}
\end{gather*}
to order any word in the algebra as wanted.
\end{proof}

\begin{Lemma}
\label{lm:basisIII}
The non-symmetric Al-Salam--Chihara polynomials form a~basis in the space $\mathcal A$ of Laurent polynomials and are
eigenfunctions of the operators~$Y$ and~$Z$:
\begin{gather*}
(Y E_{-n})[z]= \frac{1}{q^n} E_{-n}[z],
\qquad
n=1,2,\dots,
\\
(Y E_{n})[z]=0,
\qquad
n=0,1,2,\dots.
\\
(Z E_{-n})[z]=0,
\qquad
n=0,1,2,\dots,
\\
(Z E_{n})[z]=- \frac{1}{a b q^n} E_{n}[z],
\qquad
n=1,2,\dots.
\end{gather*}
\end{Lemma}

\begin{proof}
Consider the following isomorphism:
\begin{gather*}
\eta(T_0,T_1,X)=(-X T_0, T_1,X)=\big(\tilde T_0,\tilde T_1,\tilde X\big),
\end{gather*}
which maps the algebra $\mathcal H_{\rm III}$ to the isomorphic algebra $\tilde{\mathcal H}_{\rm III}$ def\/ined by the
generators $\tilde T_0$,~$\tilde T_1$,~$\tilde X$ and relations
\begin{gather*}
\big(\tilde T_1+ab\big)\big(\tilde T_1+1\big)=0,
\qquad
\tilde T_0^2+\tilde T_0=0,
\\
\big(\tilde T_1 \tilde X +a\big)\big(\tilde T_1 \tilde X +b\big)=0,
\qquad
q \tilde T_0 \tilde X^{-1}= \tilde X \big(\tilde T_0+1\big).
\end{gather*}
Note that the algebra $\tilde{\mathcal H}_{\rm III}$ is obtained by taking the limit of $c\to 0$ of the algebra
$\mathcal H_{\rm V}$, so that the proof of this lemma is based on the fact that the action of the new~$Y$ and~$Z$ def\/ined
by~\eqref{eq:new-YZ} is obtained by taking the limit of $c\to 0$ of the corresponding action of the old~$Y$ and~$Z$
def\/ined in Section~\ref{se:PV}.
\end{proof}

\begin{proof}[Proof of Theorem~\ref{th:repPIII}.] Similarly to the proof of Lemma~\ref{lm:basisIII}, we can use the
isomorphism~$\eta$ to prove this theorem by taking the limit $c\to 0$ of the proof of Theorem~\ref{th:repPV}~-- note that
this limit none of the coef\/f\/icients in the relations~\eqref{eq-ej-z} becomes zero, thus making this limit rather
straight-forward.
\end{proof}

\section[Non-symmetric continuous (big) $q$-Hermite polynomials and basic representations of ($\mathcal H_{\rm III}
^{D_7}$) $\mathcal H_{\rm III}^{D_8}$]{Non-symmetric continuous (big) $\boldsymbol{q}$-Hermite polynomials\\
and basic representations of ($\boldsymbol{\mathcal H_{\rm III}^{D_7}}$) $\boldsymbol{\mathcal H_{\rm III}^{D_8}}$}

In this section we give all def\/initions and proof for the symmetric continuous big $q$-Hermite polynomials and the
algebra $\mathcal H_{\rm III}^{D_7}$.
By taking the simple limit $a\to 0$, all proofs remain valid for the $\mathcal H_{\rm III}^{D_8}$ algebra and the
continuous $q$-Hermite polynomials.

The continuous big $q$-Hermite polynomials are the following:
\begin{gather*}
H_n(z;a):=z^n \, {}_2\phi_0\left(
\begin{matrix}
q^{-n},az
\\
-
\end{matrix}
;q,q^n z^{-2} \right),
\end{gather*}
and can be obtained from the Al-Salam--Chihara polynomials as limits when $b\to 0$.
By taking the limit $b\to 0$ of the non-symmetric continuous dual Al-Salam--Chihara we obtain the following:

\begin{Definition}
Let
\begin{gather*}
Q_n^\dagger (z;a):= q^{\frac{n-1}{2}} z  H_{n-1}\big(q^{-\frac{1}{2}} z;q^{\frac{1}{2}}a\big),
\end{gather*}
the non-symmetric continuous big $q$-Hermite polynomials are def\/ined as follows:
\begin{gather*}
E_{-n}[z]:= H_n(z;a)-Q_n^\dagger (z;a),
\qquad
n=1,2,\dots,
\\
E_{n}[z]:=q^n H_n(z;a)+\big(1-q^n\big) Q_n^\dagger (z;a),
\qquad
n=1,2,\dots,
\\
E_{0}[z]:= 1.
\end{gather*}
\end{Definition}

Similarly, the non-symmetric continuous $q$-Hermite polynomials are def\/ined by taking the limit of the non-symmetric
continuous big $q$-Hermite polynomials as $a\to 0$.

\begin{Theorem}
\label{th:repD7}
For $q,a\neq 0$, $q^m\neq 1$ $(m=1,2,\dots)$, the algebra $\mathcal H_{\rm III}^{D_7}$ has a~faithful representation on
the space $\mathcal A$ of Laurent polynomials $f[z]$ as follows:
\begin{gather}
\label{t0PIIID7}
(T_0 f)[z]:=-\frac{z}{q-z^2} \big(f[z] - f\big[q z^{-1}\big]\big),
\\
\label{t1PIIID7}
(T_1 f)[z]:=\frac{a z-1}{1-z^2} \big(f[z] - f\big[z^{-1}\big]\big),
\\
\label{XPIIID7}
(X f)[z]:=z f[z].
\end{gather}
By taking the above representation for $a=0$ $($still assuming $q\neq 0$, $q^m\neq 1$ for $m=1,2,\dots)$, we obtain
a~faithful representation of the algebra $\mathcal H_{\rm III}^{D_8}$.
\end{Theorem}

{\sloppy To prove that the operators def\/ined by~\eqref{t0PIIID7}--\eqref{XPIIID7} satisfy the
relations~\eqref{sahiD71}--\eqref{sahiD74} is a~straightforward computation.
To prove faithfulness, we can't use an equivalent representation \`a la Bernstein--Zelevinsky as there isn't one.
We proceed by proving the following two lemmata:

}

\begin{Lemma}
The algebras $\mathcal H_{\rm III}^{D_7}$ and $\mathcal H_{\rm III}^{D_8}$ are spanned by the elements
\begin{gather*}
X^k (T_0 T_1)^l,
\qquad
X^k (T_0 T_1)^lT_0,
\qquad
X^k (T_1 T_0)^l,
\qquad
X^k (T_1 T_0)^lT_1
\qquad
\hbox{for}
\quad
k\in\mathbb Z, l\in\mathbb N.
\end{gather*}
\end{Lemma}

\begin{proof}
We give the proof for the algebra $\mathcal H_{\rm III}^{D_7}$ only, as the limit $a\to 0$ in this proof is
a~straightforward substitution of~$a$ by~$0$.

Let us consider all possible words in the algebra $\mathcal H_{\rm III}^{D_7}$ and order them by using
relations~\eqref{sahiD73} and~\eqref{sahiD74} in such a~way that all powers of~$X$ are on the left.
Thanks to~\eqref{sahiD71} and~\eqref{sahiD72} the generators $T_0$ and $T_1$ may only appear with powers $1$ or~$0$.
We then are the following possible words:
\begin{gather*}
X^k (T_0 T_1)^l,
\qquad
X^k (T_0 T_1)^lT_0,
\qquad
X^k (T_1 T_0)^l,
\qquad
X^k (T_1 T_0)^lT_1
\qquad
\hbox{for}
\quad
k\in\mathbb Z, l\in\mathbb N,
\end{gather*}
as we wanted to prove.
\end{proof}

\begin{Lemma}
The non-symmetric big $q$-Hermite polynomials form a~basis in the space $\mathcal A$ of Laurent polynomials and the
operators $T_0$ and $T_1$ act on them as follows:
\begin{gather}
(T_0 E_{j}) [z]=0,
\label{eq:T0nsp}
\\
(T_0 E_{-j}) [z]=-\frac{1}{q^j} E_{j-1}[z],
\label{eq:T0nsm}
\\
(T_1 E_{j}) [z]=\big(1-q^j\big) E_{-j}[z],
\label{eq:T1nsp}
\\
(T_1 E_{-j}) [z]=- E_{-j}[z].
\label{eq:T1nsm}
\end{gather}
\end{Lemma}

\begin{proof}
By using the def\/inition of the $q$-hypergeometric series $ {}_2\phi_0$ it is easy to prove that the terms with the
highest powers in~$z$ and $\frac{1}{z}$ in $E_{-n}$ and $E_n$ have the following form
\begin{gather}
\label{andamentimD7}
E_{-n}[z]= z^{-n}+\dots - az^{n-1},
\qquad
n=1,2,\dots,
\\
\label{andamentipD7}
E_{n}[z]=z^n+\dots + q^n z^{-n},
\qquad
n=1,2,\dots.
\end{gather}
Using these relations it is straightforward to prove that the non-symmetric big $q$-Hermite polynomials form a~basis in
$\mathcal A$.

To prove~\eqref{eq:T0nsp}--\eqref{eq:T1nsm} we use the recurrence relation of the
big $q$-Hermite polynomials combined with the forward shift relation.
\end{proof}

\begin{proof}
[Proof of Theorem~\ref{th:repD7}.] To prove faithfulness we f\/irst look at how the operators $X^k (T_0 T_1)^l$, $ X^k
(T_0 T_1)^lT_0$, $X^k (T_1 T_0)^l$ and $X^k (T_1 T_0)^lT_1$ act on the non-symmetric big $q$-Hermite polynomials for
every $k\in\mathbb Z$, $l\in\mathbb N$.
To this aim, using~\eqref{eq:T0nsp}--\eqref{eq:T1nsp} one can prove the following relations:
\begin{gather*}
\big(X^k(T_0 T_1)^l E_{j}\big) [z]=-\frac{1-q^j}{q^j}\big(X^k(T_0 T_1)^{l-1} E_{j-1}\big) [z],
\qquad
\forall\, j>0,
\\
\big(X^k(T_0 T_1)^l E_{-j}\big)[z]=\frac{1}{q^j}\big(X^k(T_0 T_1)^{l-1} E_{j-1}\big) [z],
\qquad
\forall\, j>0,
\\
\big(X^k(T_0 T_1)^l T_0E_{-j}\big) [z]=\frac{1}{q^j}\frac{1-q^{j-1}}{q^{j-1}}\big(X^k(T_0 T_1)^{l-1} E_{j-2}\big) [z],
\qquad
\forall\, j>1,
\\
\big(X^k(T_1 T_0)^l E_{-j}\big) [z]=-\frac{1-q^{j-1}}{q^{j}}\big(X^k(T_1 T_0)^{l-1} E_{-j+1}\big) [z],
\qquad
\forall\, j>0,
\\
\big(X^k(T_1 T_0)^l T_1E_{j}\big) [z]=-\frac{1-q^{j}}{q^{j}}(1-q^{j-1}) \big(X^k(T_1 T_0)^{l-1} E_{-j+1}\big) [z],
\qquad
\forall\, j>0,
\\
\big(X^k(T_1 T_0)^l T_1E_{-j}\big) [z]=\frac{1-q^{j-1}}{q^{j}}\big(X^k(T_1 T_0)^{l-1} E_{-j+1}\big) [z],
\qquad
\forall\, j>0.
\end{gather*}
By iteration it is straight-forward to obtain:
\begin{gather*}
\big(X^k(T_0 T_1)^l E_{j}\big) [z]=(-1)^l\frac{\left(q^{j-l+1};q\right)_l}{q^{\frac{l(1+2j-l)}{2}}} \big(X^k E_{j-l}\big) [z],
\qquad
\forall\, j\geq l,
\\
\big(X^k(T_0 T_1)^l E_{-j}\big)[z]=(-1)^{l-1}\frac{\left(q^{j-l+1};q\right)_{l-1}}{q^{\frac{l(1+2j-l)}{2}}} \big(X^k E_{j-l}\big) [z],
\qquad
\forall\, j\geq l,
\\
(X^k(T_0 T_1)^l T_0E_{-j}) [z]=(-1)^{l-1}\frac{\left(q^{j-l};q\right)_{l}}{q^{\frac{(l+1)(2j-l)}{2}}} (X^k E_{j-l-1})[z],
\qquad
\forall\, j> l,
\\
\big(X^k(T_1 T_0)^l E_{-j}\big) [z]=(-1)^{l}\frac{\left(q^{j-l};q\right)_{l}}{q^{\frac{l(1+2j-l)}{2}}} \big(X^k E_{-j+l}\big) [z],
\qquad
\forall\, j\geq l,
\\
\big(X^k(T_1 T_0)^l T_1E_{j}\big) [z]=(-1)^{l}\frac{\left(q^{j-l};q\right)_{l+1}}{q^{\frac{l(1+2j-l)}{2}}} \big(X^k E_{-j+l}\big) [z],
\qquad
\forall\, j\geq l,
\\
\big(X^k(T_1 T_0)^l T_1E_{-j}\big) [z]=(-1)^{l-1}\frac{\left(q^{j-l};q\right)_{l}}{q^{\frac{l(1+2j-l)}{2}}} \big(X^k E_{-j+l}\big) [z],
\qquad
\forall\, j\geq l.
\end{gather*}
Combining these with~\eqref{andamentimD7} and~\eqref{andamentipD7}, we obtain the following estimates $\forall\, j>l$:
\begin{gather*}
\big(X^k(T_0 T_1)^l E_{j}\big) [z]=(-1)^l\frac{\big(q^{j-l+1};q\big)_l}{q^{\frac{l(1+2j-l)}{2}}} \big(z^{k+j-l}+\dots +q^{j-l}z^{k-j+l}\big),
\\
\big(X^k(T_0 T_1)^l E_{-j}\big)[z]=(-1)^{l-1}\frac{\big(q^{j-l+1};q\big)_{l-1}}{q^{\frac{l(1+2j-l)}{2}}} \big(z^{k+j-l}+\dots+q^{j-l} z^{k-j+l}\big),
\\
\big(X^k(T_0 T_1)^l T_0E_{-j}\big) [z]=(-1)^{l-1}\frac{\big(q^{j-l};q\big)_{l}}{q^{\frac{(l+1)(2j-l)}{2}}} \big(z^{k+j-l-1}+\dots+q^{j-l} z^{k-j+l+1}\big),
\\
(X^k(T_1 T_0)^l E_{-j}) [z]=(-1)^{l}\frac{\big(q^{j-l};q\big)_{l}}{q^{\frac{l(1+2j-l)}{2}}} (z^{k-j+l}+\dots-az^{k+j-l-1}),
\\
\big(X^k(T_1 T_0)^l T_1E_{j}\big) [z]=(-1)^{l}\frac{\big(q^{j-l};q\big)_{l+1}}{q^{\frac{l(1+2j-l)}{2}}}\big(z^{k-j+l}+\dots-az^{k+j-l-1}\big),
\\
\big(X^k(T_1 T_0)^l T_1E_{-j}\big) [z]=(-1)^{l-1}\frac{\big(q^{j-l};q\big)_{l}}{q^{\frac{l(1+2j-l)}{2}}} \big(z^{k-j+l}+\dots-az^{k+j-l-1}\big).
\end{gather*}
Now assume by contradiction that a~linear combination acts as zero operator in our representation, let us write such
linear combination as:
\begin{gather*}
\sum\limits_{k,l} a_{k,l} X^k (T_0 T_1)^l+\sum\limits_{k,l}b_{k,l}X^k(T_0 T_1)^lT_0+\sum\limits_{k,l}c_{k,l}X^k (T_1T_0)^l
+ \sum\limits_{k,l} d_{k,l}X^k (T_1 T_0)^lT_1.
\end{gather*}
Take the minimum value $k_0$ of~$k$ such that at least one coef\/f\/icient $a_{k_0,l}$, $b_{k_0,l}$, $c_{k_0,l}$, $d_{k_0,l}$ is nonzero.
Acting on $E_j[z]$, for all $j>l$, and collecting the terms with the minimum possible power of~$z$, which is
$z^{k_0-j+l}$, we obtain the equation:
\begin{gather*}
\sum\limits_{l} a_{k_0,l} (-1)^l\frac{\big(q^{j-l+1};q\big)_l}{q^{\frac{l(1+2j-l)}{2}}} q^{j-l} + \sum\limits_{l}
d_{k_0,l}(-1)^{l}\frac{\big(q^{j-l};q\big)_{l+1}}{q^{\frac{l(1+2j-l)}{2}}}=0,
\qquad
\forall\, j>l.
\end{gather*}
It is easy to prove that for generic values of~$a$, this is an inf\/inite set of linearly independent equations, therefore
the only possible solution is the trivial one, i.e.~$a_{k_0,l}=0$, $d_{k_0,l}=0$ for all values of~$l$.

By acting on $E_{j}[z]$, we can prove in a~similar way that $b_{k_0,l}=0$, $c_{k_0,l}=0$ for all values of~$l$, therefore
obtaining a~contradiction.

To prove the same for the algebra $\mathcal H_{\rm III}^{D_8}$ we observe that the def\/ining
relations~\eqref{sahiPIIID82}--\eqref{sahiPIIID85} are a~specialisation of the def\/ining
relations~\eqref{sahiD71}--\eqref{sahiD74} of the algebra $\mathcal H_{\rm III}^{D_7}$ for $a=0$.
All results hold true when $a\to 0$.
Indeed even if the polynomials~$E_{n}$ loose the terms of order~$z^{n-1}$, these don't enter in the above reasoning.
This concludes the proof of our theorem.
\end{proof}

\subsection*{Acknowledgements}
The author is specially grateful to T.~Koornwinder and J.~Stokman for interesting discussions on this subject and the
referees for pointing out references, typos and giving suggestions on how to improve the presentation of the main
results.

\pdfbookmark[1]{References}{ref}
\LastPageEnding

\end{document}